\documentclass[preprint,12pt,authoryear]{elsarticle}

\usepackage{amssymb}

\usepackage{amsmath}
\usepackage{changepage}
\usepackage{physics}
\usepackage{subcaption} 
\usepackage{graphicx} 
\usepackage{doi}
\usepackage{float}
\usepackage{booktabs} 
\usepackage{hyperref}
\usepackage{natbib}
\journal{Nuclear Physics B}
\usepackage{grffile}  % Permette spazi e caratteri speciali nei nomi

\begin{document}

\thispagestyle{plain}

\begin{center}
    \LARGE{\textbf{Quantum-Enhanced Spectral Solution of the Poisson Equation}}
\end{center}

\begin{center}
\vspace{4pt}
\large
    G. Intoccia\textsuperscript{1}, U. Chirico \textsuperscript{1,2},  G. Pepe\textsuperscript{1}, S. Cuomo \textsuperscript{1} \\
    
\vspace{0.4cm}
\small
   \textsuperscript{1} University of Naples Federico  \textsuperscript{2} Quantum2pi S.r.l. 

\end{center}

\begin{small}
\begin{center}
\vspace{9pt}
\textbf{Abstract}    
\end{center}

\begin{adjustwidth}{20pt}{20pt}
\small \noindent We present a hybrid numerical-quantum method for solving the Poisson equation under homogeneous Dirichlet boundary conditions, leveraging the Quantum Fourier Transform (QFT) to enhance computational efficiency and reduce time and space complexity. This approach bypasses the integration-heavy calculations of classical methods, which have to deal with high computational costs for large number of points.

The proposed method estimates the coefficients of the series expansion of the solution directly within the quantum framework. Numerical experiments validate its effectiveness and reveal significant improvements in terms of time and space complexity and solution accuracy, demonstrating the capability of quantum-assisted techniques to contribute in solving partial differential equations (PDEs).

Despite the inherent challenges of quantum implementation, the present work serves as a starting point for future researches aimed at refining and expanding quantum numerical methods.
\end{adjustwidth}

\end{small}

\vspace{11pt}
\section{Introduction}\label{Sec:Introduction}
The Poisson equation arises in various physical contexts, such as electrostatics, fluid dynamics, and heat conduction. Solving this equation with Dirichlet boundary conditions is a classical problem that has been extensively studied using various numerical and analytical techniques.

In recent years, the growing interest in quantum computing has opened new perspectives for addressing complex mathematical problems, including partial differential equations (PDEs). Classical numerical methods, such as finite elements or spectral methods, are powerful but often suffer from computational limitations, especially for high-dimensional problems or geometrically complex domains. These limitations have motivated the exploration of hybrid approaches that use the strengths of quantum computation.
\\
The proposed solution to the Poisson equation combines classical numerical techniques with quantum computation to estimate the coefficients from the series expansion of the solution.
We adopt a probabilistic approach to the spectral series coefficients; specifically, by measuring the quantum state after the application of the Quantum Fourier Transform (QFT), we obtain the probability of observing frequencies that are present in the original state. In other words, the measurement yields a probability distribution corresponding to the  modulus-squared of the spectral coefficients.\\
By using the QFT, the method avoids the need for classical techniques that require the evaluation of multiple integrals, which are computationally expensive. This quantum-enhanced approach reduces computational complexity and enables more efficient solutions.

\section{Quantum Method Implementation for 2D Poisson Problems with homogeneous Dirichlet boundary conditions} 

In this article, we focus on solving the Poisson equation in two dimensions, specifically in a rectangular domain:

\begin{equation}
\begin{cases} 
\Delta u = f & \text{in } \Omega \\
u = 0 & \text{on } \partial \Omega
\end{cases}
\end{equation}

where \( \Omega = [0,L_x] \times [0, L_y] \) and \( \Delta  = \frac{\partial^2 }{\partial x^2} + \frac{\partial^2 }{\partial y^2} \).\\
The solution \( u(x, y) \) is expressed as a linear combination of the eigenfunctions \( \phi_{k_xk_y}(x, y) \), given by:

\begin{equation}
    \phi_{k_xk_y}(x,y) = \sin\left(\frac{k_x \pi x}{L_x}\right) \sin\left(\frac{k_y \pi y}{L_y}\right)
\end{equation}

where \( k_x,k_y \in \mathbb{N} \).

The general solution is:
\begin{equation}
    u(x,y) = \sum_{k_x=0}^\infty \sum_{k_y=0}^\infty  a_{k_x k_y} \sin\left(\frac{(k_x+1) \pi x}{L_x}\right) \sin\left(\frac{(k_y+1) \pi y}{L_y}\right)
\end{equation}

where the coefficients \( a_{k_xk_y} \) are determined as:

\begin{equation}
    a_{k_xk_y} = 4 \cdot \frac{\int_\Omega f(x,y) \sin\left(\frac{k_x \pi x}{L_x}\right) \sin\left(\frac{k_y \pi y}{L_y}\right) dx dy}{L_x L_y}
\end{equation}

and \( f(x,y) \) is the source function in the Poisson equation.

We construct the solution to this equation using a hybrid quantum-classical method. First, we discretize the domain:
\begin{align*}
    x_i &= i \frac{L_x}{2^n}, \quad i = 0, \ldots, N  \\
    y_j &= j \frac{L_y}{2^m}, \quad j = 0, \ldots, M 
\end{align*}

\( n \) and \( m \) are the number of qubits used to represent the \( x \)- and \( y \)-dimensions of a 2D spatial domain, respectively, with $N=2^n$ ,  $M=2^m$. \\ The corresponding Hilbert space is given by:

\begin{equation}
    \mathcal{H}_x \otimes \mathcal{H}_y \simeq \mathbb{C}^{2^{n+m}}
\end{equation}

To map the function $f(x,y)$ into a quantum state, we first normalize it to ensure compatibility with the quantum state initialization:
\begin{equation}
    \psi(x,y)=\frac{f(x,y)}{||f||_2}
\end{equation}
where $||f||_2$ is the $L^2([0,L_x] \times [0, L_y])$ norm.
Then we construct a vector by flattening the values of $\psi(x,y)$:
\begin{equation*}
    \mathbf{\psi(x,y)} = \Big(\psi(x_0,y_0), \dots,\psi(x_0,y_{N}), \psi(x_1,y_0), \psi(x_1,y_1),\dots,\psi(x_1,y_{M}), \dots \Big)
\end{equation*}
Renaming the components for clarity:
\begin{equation*}
    \mathbf{\psi(x,y)} = \Big(\psi_{00},\dots,\psi_{0N}, \psi_{10},\dots,\psi_{1M}, \dots,\psi_{ij},\dots \Big)
\end{equation*}

The quantum state $\ket{\psi}\in \mathcal{H}_x \otimes \mathcal{H}_y$  is defined as:

\begin{equation*}
    \ket{\psi} = \sum_{i=0}^{N - 1} \sum_{j=0}^{M- 1} \psi_{ij} \, \ket{i} \otimes \ket{j}
\end{equation*}

where:

\begin{itemize}
    
    \item $ \ket{i}$ is a computational basis state of the \( n \)-qubit register representing the \( x \)-coordinate index, \( i \in \{0, \dots, 2^n - 1\} \);
    \item \( \ket{j} \) is a computational basis state of the \( m \)-qubit register representing the \( y \)-coordinate index, \( j \in \{0, \dots, 2^m - 1\} \);
    \item \( \psi_{ij} \in \mathbb{C} \) the components of $\mathbf{\psi(x,y)}$ that represents the amplitudes of the state.
\end{itemize}

\subsection{Quantum Fourier Transform}
The Quantum Fourier Transform (QFT) is the quantum analogue of the Discrete Fourier Transform (DFT). For an \(n\)-qubit state, the QFT is defined as:
\begin{equation*}
    \text{QFT}|k\rangle = \frac{1}{\sqrt{2^n}} \sum_{j=0}^{2^n - 1} e^{2\pi i k j / 2^n} |j\rangle
\end{equation*}
where \(|k\rangle\) is the computational basis state.

The QFT circuit applies a series of Hadamard gates and controlled phase rotations to create this transformation efficiently.\\
In our approach we apply the QFT first to the $n$ qubits, and then to the $m$ qubits.\\
The resulting state is:
\begin{equation}
    |\tilde{\psi}\rangle = (\text{QFT}_n \otimes \text{QFT}_m) |\psi\rangle
\end{equation}

The new amplitudes \( \tilde{\psi}_{kl} \) are given by the 2D discrete Fourier transform:

\begin{equation}
    \tilde{\psi}_{kl} = \frac{1}{\sqrt{2^{n+m}}} \sum_{w=0}^{2^n - 1} \sum_{j=0}^{2^m - 1} \psi_{ij} \, e^{2\pi i \frac{wk}{2^n}} \, e^{2\pi i \frac{jl}{2^m}}
\end{equation}

and the state becomes:

\begin{equation}
    |\tilde{\psi}\rangle = \sum_{k=0}^{2^n - 1} \sum_{l=0}^{2^m - 1} \tilde{\psi}_{kl} \, |k\rangle \otimes |l\rangle
\end{equation}

\subsection{Measurement and Spectral Estimation}
Upon measuring the quantum state after the 2D QFT, the probability of obtaining the binary outcome $\ket{kl}$ is:

\begin{equation}
    \mathbb{P}(k, l) = |\tilde{\psi}_{kl}|^2
\end{equation}

which represents the frequencies distribution of the original state amplitudes \( \psi_{ij} \).

Let \( C_{kl} \) be the number of occurrences of the state $\ket{kl}$, corresponding to the integer indices $(k,l)$ obtained from binary measurement outcomes. The empirical probability is then estimated as

\begin{equation}
    \mathbb{P}(k, l) \approx \frac{C_{kl}}{S}
\end{equation}

where $S$ are the shots. The magnitude of the frequencies amplitude is approximated by

\begin{equation}
    |\tilde{\psi}_{kl}| \approx \sqrt{\frac{C_{kl}}{S}}
\end{equation}

We gather the counts in a matrix \( a_{kl} \)

\begin{equation}
    a_{kl} \propto \frac{\sqrt{C_{kl}}}{\sqrt{\sum_{k,l} C_{kl}}} 
\end{equation}

This normalization guarantees that the reconstructed state satisfies

\begin{equation}
    \sum_{k,l} |a_{kl}|^2 = 1
\end{equation}

We define the eigenvalues for the Laplacian on the rectangular domain as
\begin{equation*}
    \lambda_x(i) = \left(\frac{\pi (i+1)}{L_x}\right)^2 \hspace{0.5cm} i=0,1,\ldots,N-1
\end{equation*}

\begin{equation*}
    \lambda_y(j) = \left(\frac{\pi (j+1)}{L_y}\right)^2 \hspace{0.5cm} j=0,1,\ldots,M-1
\end{equation*}

The total eigenvalue is given by
\begin{equation}
    \lambda_{ij} = \lambda_x(i) + \lambda_y(j)
\end{equation}
and the final coefficients are given by
\begin{equation}
    b_{ij}=\frac{a_{ij}}{\lambda_{ij}}
\end{equation}
The eigenfunctions (modes) for indices \(i,j\) are given by the sine basis:
\begin{equation}
    \phi_{ij}(x,y) = \sin\left(\frac{\pi (i+1) x}{L_x}\right) \sin\left(\frac{\pi (j+1) y}{L_y}\right)
\end{equation}

Then, the solution $u(x,y)$ is reconstructed as a linear combination of these modes weighted by the coefficients $a_{ij}$ divided by the eigenvalues:
\begin{equation}
    u(x,y) = \sum_{i=0}^{N - 1} \sum_{j=0}^{M - 1} \frac{\operatorname{Re}(a_{ij})}{\lambda_{ij}} \phi_{ij}(x,y)
\end{equation}

Since the coefficients $a_{ij}$ can be complex due to the quantum state, we take their real part in the final reconstruction.

\section{Experimental Results}

In this section, we present the experimental results obtained using the quantum code described previously and compare them with the results generated using the classical Fourier series method,,so with the direct integration of the coefficients by \texttt{dblquad} from python. The \texttt{AerSimulator()} from Qiskit is used to simulate the behavior of the quantum circuit.\\
The domain, for all the cases that we are going to consider, is fixed to $[0,1]^2$.\\
We first consider a purely sinusoidal source function:
\begin{equation*}
    f(x,y)=\sin(k_1 \pi x) \sin(k_2 \pi y)\         \
\end{equation*}

To ensure correct harmonic retrieval, a shifting procedure was applied:

\begin{equation*}
     k_1 = m_1 + \text{shift}(m_1) \hspace{0.7cm}  k_2 = m_2 + \text{shift}(m_2) \hspace{0.7cm} m_1,m_2=1,2,3,\dots 
\end{equation*}

where:

\begin{equation}
\text{shift}(\text{m}) = 
\begin{cases} 
0 & \text{if } \text{m} \leq 2 \\
\text{m} - 2 & \text{otherwise}
\end{cases}
\end{equation}
After obtaining the coefficients from Eq.(13), a correction was applied to address the loss of relative phase information after the measurements. The problem was particularly evident for specific modes. By applying a phase and spectral correction, the accurate reconstruction of the solution was recovered.
It is also worth noticing that all experimental validations were carried out using a limited number of qubits, ranging from 10 to 16 in total.
\begin{equation}
B_{mn} = a_{mn} \cdot \underbrace{\exp\left(-i\pi\frac{m + n}{2}\right)}_{\text{Phase}} \cdot \underbrace{\frac{(mn)^2}{(m + n)^3}}_{\text{Spectral correction}}
\end{equation}

The results are shown in Figure 1.

\begin{figure}[H]
  \centering
  \begin{subfigure}[b]{0.48\textwidth}
    \includegraphics[width=\textwidth]{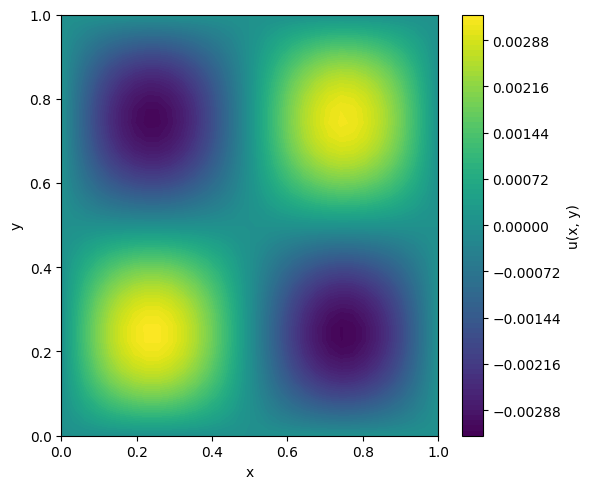}
    \caption{}
    \label{fig:sub1}
  \end{subfigure}
  \hfill
  \begin{subfigure}[b]{0.48\textwidth}
    \includegraphics[width=\textwidth]{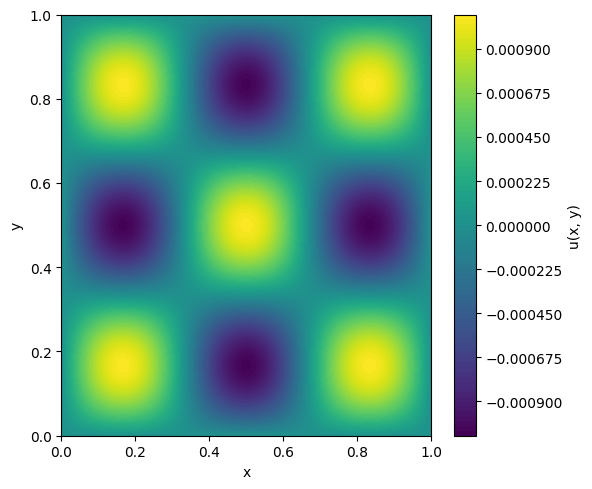}
    \caption{}
    \label{fig:sub2}
  \end{subfigure}
  \caption{Results for (a) Left: $(k_1,k_2)=(2,2)$ and (b) Right: $(k_1,k_2)=(3,3)$.}
  \label{fig:global}
\end{figure}

The mean squared errors (MSE) compared to the classical solution is $mse = 3.20\times 10^{-10}$ and $mse = 3.20\times 10^{-8}$, respectively.

Let's consider the following function
\begin{equation*}
    f(x,y)=x(1-x)y(1-y)
\end{equation*}

Using the same correction from the Eq.(20) we obtain the following result.
\begin{figure}[H]
  \centering
    \includegraphics[width=0.60\textwidth]{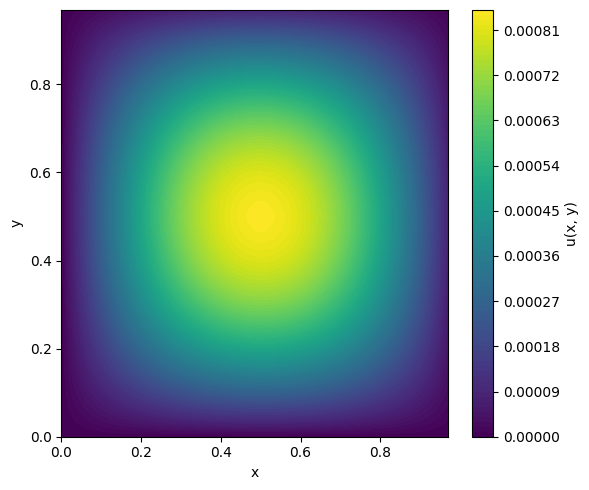}
    \caption{}
    \label{fig:sub1}

  \label{fig:global}
\end{figure}

The mean squared error is $mse = 1.32\times 10^{-9}$.

Next, we consider a spatially anisotropic source function:
\begin{equation}
f(x,y) = \sin(k_1\pi x)\sin(k_2\pi y) \cdot (x^2 + 2xy + 3y^2 - x + 4y + 5)
\end{equation}

Thus, the polynomial introduces spatial anisotropy while satisfying homogeneous Dirichlet boundary conditions, generating spectral coefficients with additional contributions compared to the purely sinusoidal case. The corrective term applied is
\begin{equation}
    B_{mn} = a_{mn} \cdot \exp\left(-i\pi(m + n)\right) \cdot \frac{(mn)^2}{(m + n)^3}
\end{equation}

\begin{figure}[H]
  \centering
  \begin{subfigure}[b]{0.48\textwidth}
    \includegraphics[width=\textwidth]{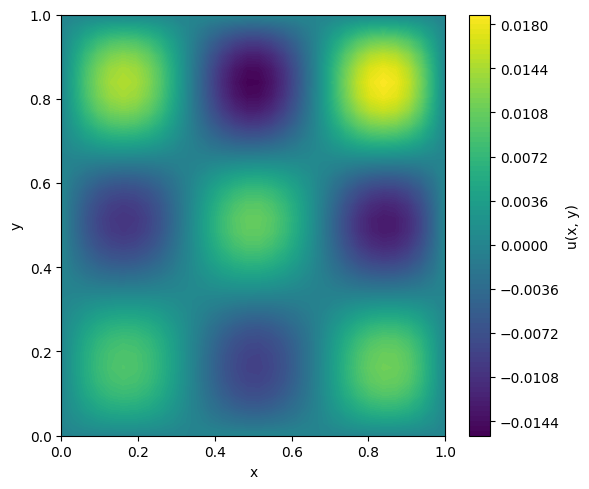}
    \caption{}
    \label{fig:sub1}
  \end{subfigure}
  \hfill
  \begin{subfigure}[b]{0.48\textwidth}
    \includegraphics[width=\textwidth]{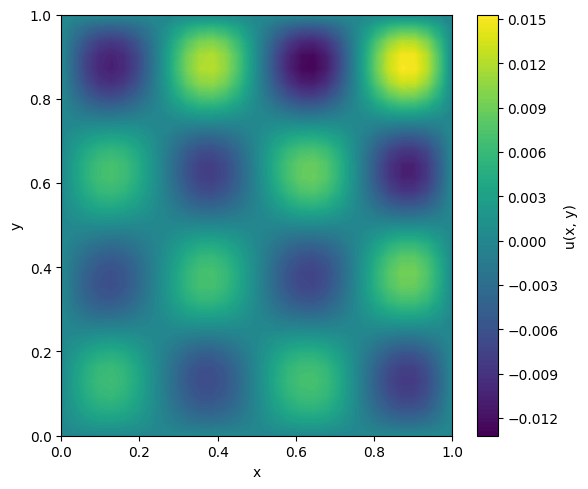}
    \caption{}
    \label{fig:sub2}
  \end{subfigure}
  \caption{Solutions for (a) Left: $(k_1,k_2)=(3,3)$ (b) Right: $(k_1,k_2)=(4,4)$.}
  \label{fig:global}
\end{figure}

The results yielded MSEs of $mse = 3.05\times 10^{-6}$ and $mse = 2.91\times 10^{-6}$, respectively.

The Gaussian source represents a critical test case for estimating the relative phases lost during the quantum computation process. Let
\begin{equation}
f(x,y) = \exp\left(-((x-x_0)^2 + (y-y_0)^2)\right) 
\end{equation}
with $(x_0,y_0)=(0.5,0.5)$.
To recover the correct solution, the coefficients must be multiplied by a phase factor. In this case, the optimal phase factor determined experimentally was found to be:
\begin{equation}
    B_{mn} = a_{mn} \cdot \exp\left(-i\pi\frac{mn}{2}\right)
\end{equation}
The solution is shown in the Figure 4.
\begin{figure}[H]
  \centering
    \includegraphics[width=0.60\textwidth]{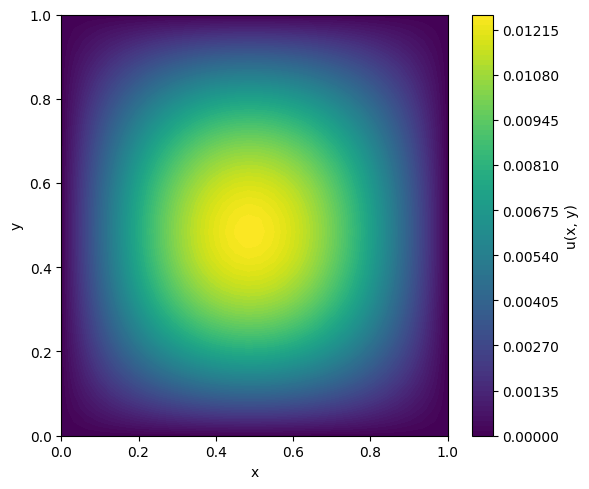}
    \caption{Solution for a Gaussian function}
\end{figure}
The resulting solution achieved $mse = 9.31\times 10^{-6}$.

Let us now consider the following product as a known term
\begin{equation}
    f(x,y) = \exp\left(-((x-x_0)^2 + (y-y_0)^2)\right)+ \sin(k_1\pi x)\sin(k_2\pi y) 
\end{equation}
In this case the optimal phase factor determined experimentally was found to be:
\begin{equation}
    B_{mn} = a_{mn} \cdot \exp\left(-i\pi mn \right)\frac{(mn)^2}{(m^2 + n^2)^\frac{3}{2}}
\end{equation}

\begin{figure}[H]
  \centering
  \begin{subfigure}[b]{0.48\textwidth}
    \includegraphics[width=\textwidth]{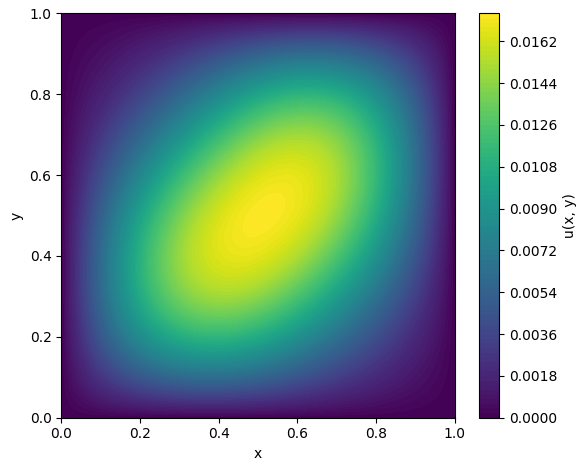}
    \caption{}
    \label{fig:sub1}
  \end{subfigure}
  \hfill
  \begin{subfigure}[b]{0.48\textwidth}
    \includegraphics[width=\textwidth]{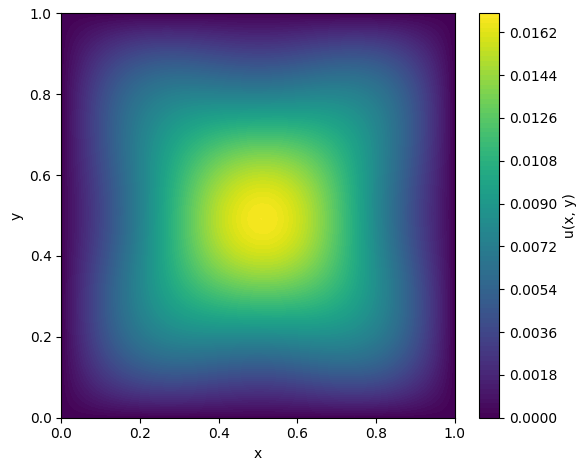}
    \caption{}
    \label{fig:sub2}
  \end{subfigure}
  \caption{Solutions for $(x_0.y_0)=(0.5,0.5)$ (a) Left: $(k_1,k_2)=(2,2)$ (b) Right: $(k_1,k_2)=(3,3)$}
  \label{fig:global}
\end{figure}

The results are shown in Figure 5 and yielded respectively $mse = 1.52\times 10^{-6}$ and $mse = 5.87\times 10^{-6}$.

\section{The advantages of using Quantum}

The determination of Fourier coefficients in a two-dimensional problem requires a detailed analysis of computational complexity. We will examine the costs associated with classical and quantum methods for calculating these coefficients.

\subsection{Computational complexity of coefficient estimation}

The direct calculation of Fourier coefficients in a discrete domain requires the execution of double sums for every combination of coefficients. Considering a grid of size $2^n \times 2^m$, the computational cost is:

\begin{equation*}
O(2^n \cdot 2^m) = O(2^{n+m})
\end{equation*}

This approach quickly becomes inefficient for large values of $n$ e $m$. To improve performance, the Fast Fourier Transform (FFT) is used, which reduces the complexity to:

\begin{equation*}
O(2^n \log(2^n) + 2^m \log(2^m))
\end{equation*}

Despite the improvement over direct calculation, the FFT approach is still limited by the exponential growth of the number of points.

The Quantum Fourier Transform (QFT) leverages the properties of superposition and entanglement to achieve a significant computational advantage. The QFT has polynomial complexity with respect to the number of qubits

\begin{equation*}
O(n^2 + m^2)
\end{equation*}

The following table summarizes the computational costs.

\begin{table}[H]
\centering
\begin{tabular}{lrrrr}
\toprule
Method & Complexity & Notes \\
\midrule
Direct Classical  & $O(2^{n+m})$ & Inefficient  \\
Classical FFT  & $O(2^n \log(2^n) + 2^m \log(2^m))$ & Improved  \\
QFT  & $O(n^2 + m^2)$ & Quasi linear \\
\bottomrule
\end{tabular}
\caption{Comparison of computational complexity.}
\end{table}

\subsection{Computational Efficiency and Memory Analysis
}

Comparisons between methods were performed on a test function.
\begin{equation*}
    f(x,y)=\sin(\pi x)\sin(\pi y)
\end{equation*}
The classical implementation employs spectral reconstruction of the solution using a standard spectral method, where coefficient computation is performed through a numerical approximation of the integral. While the results of this section are excellent for simple approximation and quadrature methods, they will prove even more effective for more advanced quadrature methods.\\
First, we observe that since the coefficients are constructed via the QFT, and thus derived from the frequencies of the original quantum state, the series can be truncated during the solution reconstruction phase
\begin{equation*}
    u(x,y) = \sum_{i=0}^{\tau_x - 1} \sum_{j=0}^{\tau_y - 1} \frac{\operatorname{Re}(a_{ij})}{\lambda_{ij}} \phi_{ij}(x,y)
\end{equation*}
where  $\tau_x$ and $\tau_y$ are powers of 2. 
This always allows us to reduce the execution time of the final phase, provided that we are careful not to truncate excessively to avoid errors. Figure 6 shows the solution for $(\tau_x,\tau_y)=(9,9)$ as a function of $n$,$m$  number of qubits.

\begin{figure}[H]
  \centering
  \begin{subfigure}[b]{0.32\textwidth}  % Larghezza ridotta per 3 figure
    \includegraphics[width=\textwidth]{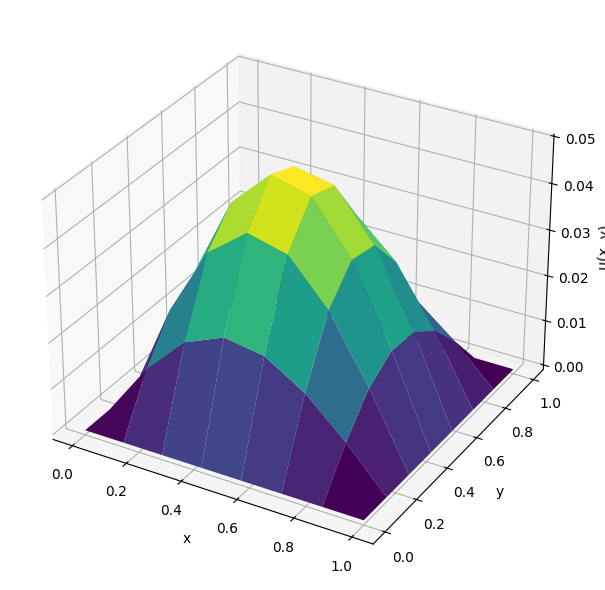}
    \caption{$(n,m)=(3,3)$}
    \label{fig:sub1}
  \end{subfigure}
  \hfill
  \begin{subfigure}[b]{0.32\textwidth}
    \includegraphics[width=\textwidth]{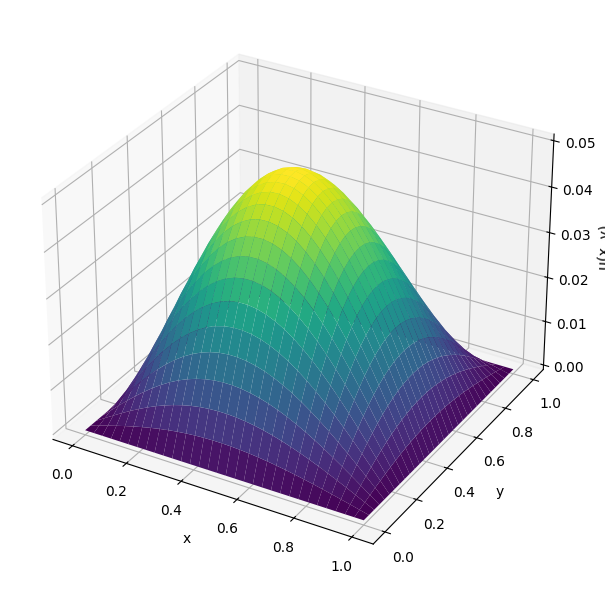}
    \caption{$(n,m)=(5,5)$}
    \label{fig:sub2}
  \end{subfigure}
  \hfill
  \begin{subfigure}[b]{0.32\textwidth}
    \includegraphics[width=\textwidth]{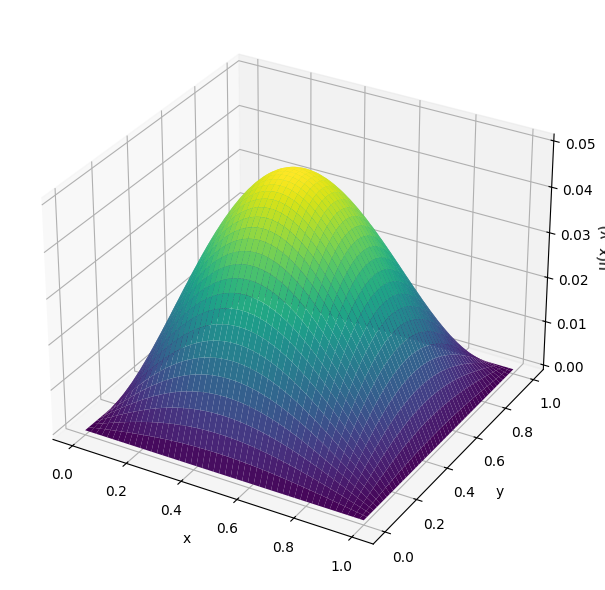}  % Sostituisci con il nome del file
    \caption{$(n,m)=(10,10)$}  % Personalizza la didascalia
    \label{fig:sub3}
  \end{subfigure}
  \caption{Truncated series solutions for $n$,$m$ number of qubits.}
  \label{fig:global}
\end{figure}

In the subsequent analysis, we will consider the untruncated series to better observe the scaling behavior.\\
The benchmark results demonstrate substantial improvements in time and space complexity when using the hybrid quantum-inspired approach compared to the classical method.
The following results refer to a $256 \times 256$ grid, which corresponds to a total of 16 qubits.

%DA RISCRIVERE
The computational performance is analyzed in Figures 7-8 (execution time, shown in both absolute and percentage terms) and Figures 9-10 (memory usage).\\
\textbf{NOTE:} in the coefficient correction phase is included the ratio between the coefficients and the eigenvalues, as indicated in Eq.(16)

\begin{figure}[H]
  \centering
    \includegraphics[width=1\textwidth]{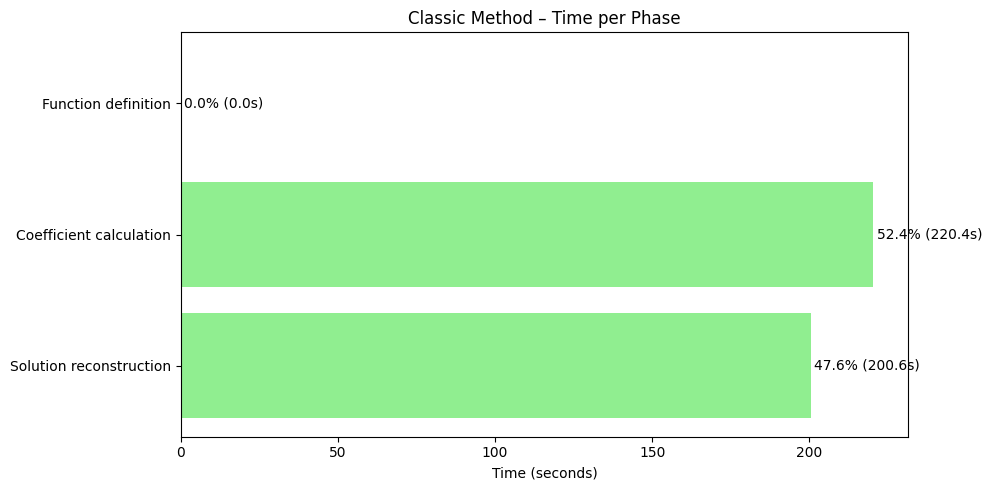}
    \caption{Execution times breakdown for classical computation phases.}
\end{figure}
\begin{figure}[H]
  \centering
    \includegraphics[width=1\textwidth]{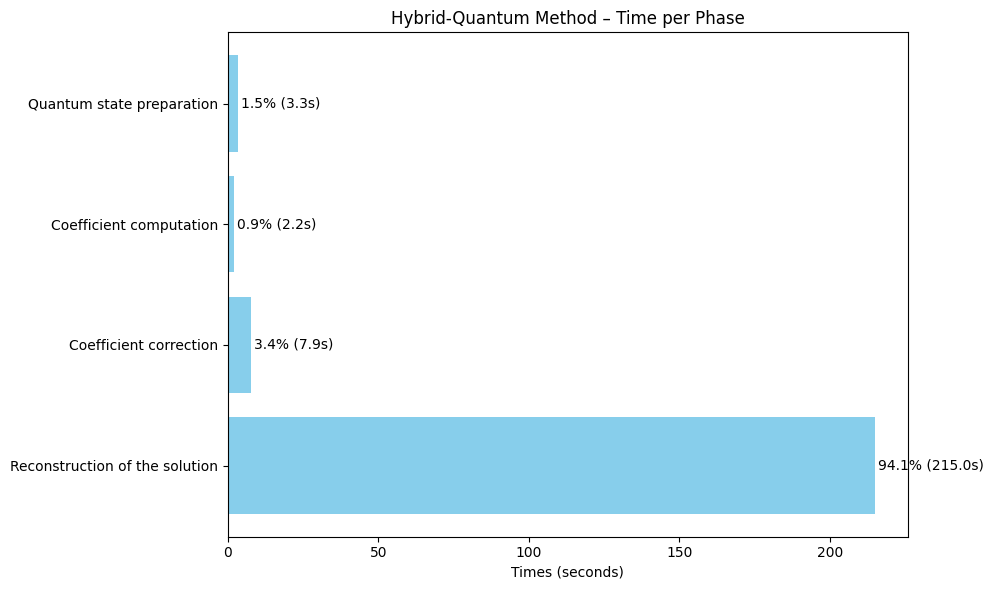}
    \caption{Execution times breakdown for hybrid-quantum computation.}
\end{figure}

%DA RISCRIVERE

\begin{figure}[H]
  \centering
    \includegraphics[width=1\textwidth]{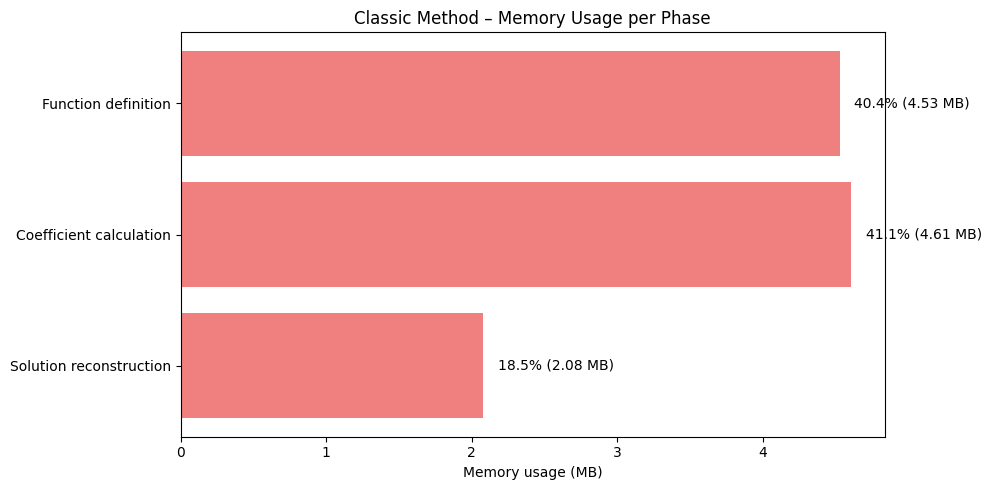}
    \caption{Memory usage distribution across classical computation phases.}
\end{figure}
\begin{figure}[H]
  \centering
    \includegraphics[width=1\textwidth]{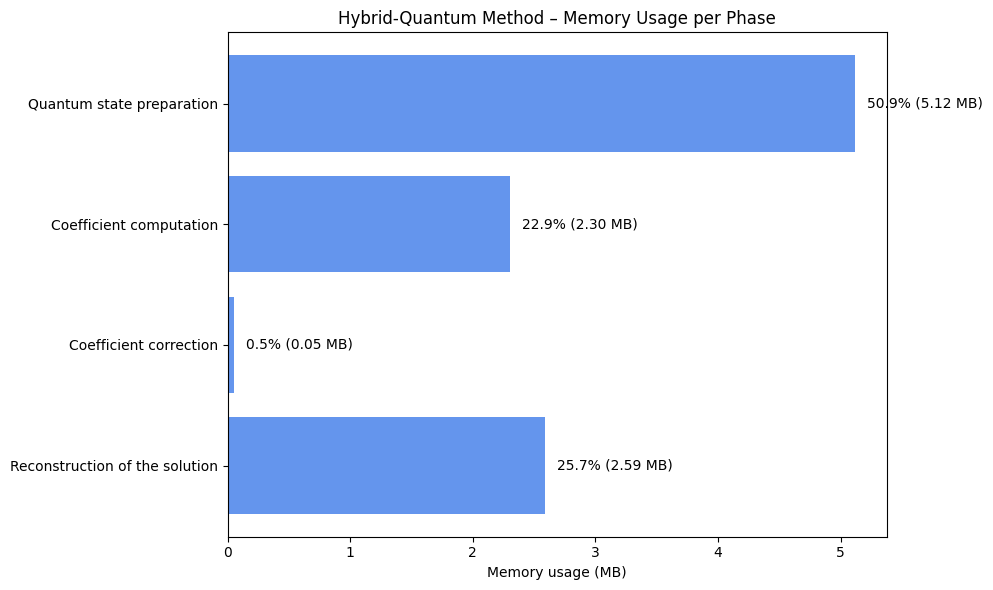}
    \caption{Memory usage distribution in hybrid-quantum computation.}
\end{figure}

The corresponding data are shown in Table 2. Regarding the hybrid-quantum case, we can note a drastic reduction of 93.9\% in the computation time of the coefficient. The reconstruction step, which remains a classical computation, becomes the dominant cost in the quantum phase, unless we apply truncation.

\begin{table}[H]
\centering
\begin{tabular}{lrrr}
\toprule
\textbf{Phase} & \textbf{Classical (s)} & \textbf{Quantum (s)} & \textbf{Change} \\
\midrule
Coefficient Calculation & 220.43 & 2.39 & 93.9\% $\downarrow$ \\
Solution Reconstruction & 200.59 & 215.04 & 7.2\% $\uparrow$ \\
State Preparation & N/A & 3.34 & -- \\
\hline
Total Time & 421.03 & 228.43 & 45.7\% $\downarrow$ \\
\bottomrule
\end{tabular}
\caption{Computational time comparison between classical and quantum approaches. The quantum method demonstrates a 45.7\% overall speedup, with particular improvement (93.9\%) in coefficient calculation. The state preparation represents additional quantum-specific overhead.}
\label{tab:performance}
\end{table}

In terms of memory usage (see Table 3), there is a slight overall reduction. While the coefficient computation phase shows a significant decrease in memory consumption, the state preparation phase requires more memory. 

\begin{table}[H]
\centering
\begin{tabular}{lrrr}
\toprule
\textbf{Phase} & \textbf{Classical (MB)} & \textbf{Quantum (MB)} & \textbf{Change} \\
\midrule
Initialization & 4.53 & 5.11 & +12.8\% $\uparrow$ \\
Coefficient Processing & 5.61 & 2.77 & -50.6\% $\downarrow$ \\
Final Phase & 2.08 & 2.59 & +24.5\% $\uparrow$ \\
\hline
Total Usage & 12.22 & 10.47 & -14.3\% $\downarrow$ \\
\bottomrule
\end{tabular}
\caption{Memory Usage Comparison Between Classical and Quantum Approaches. The quantum method shows an overall 14.3\% reduction in total memory consumption despite increased requirements in some phases.}
\label{tab:memory}
\end{table}

Figure 11 shows a comparison between the classical coefficient calculation and the quantum coefficient computation, though this pertains only to the application of the QFT. The plot shows the variation of the time execution as a function of the chosen number of points.

\begin{figure}[H]
  \centering
    \includegraphics[width=1\textwidth]{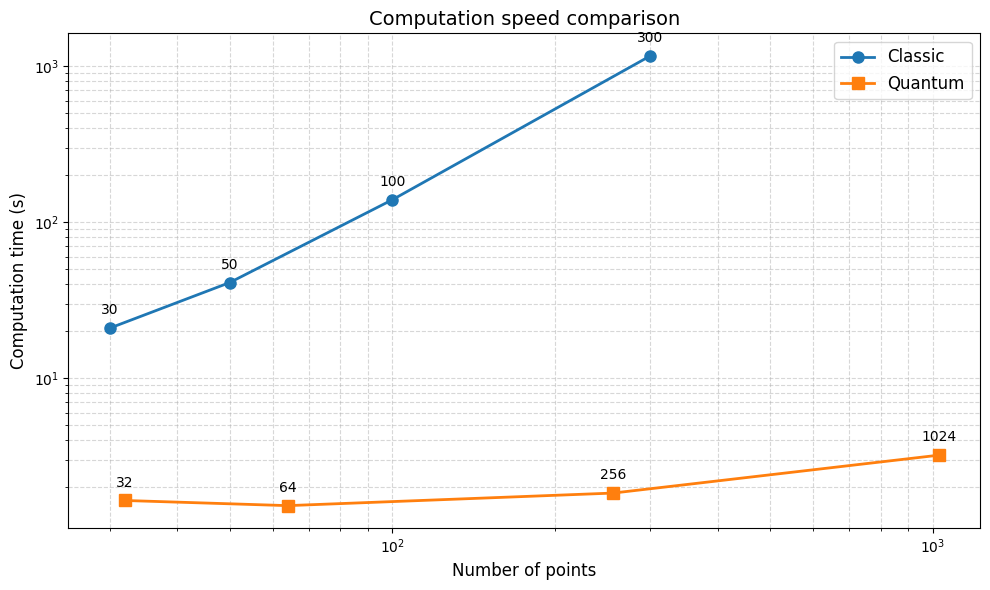}
    \caption{Runtime Classical vs. QFT-based quantum coefficient computation.}
\end{figure}
However, for a fair comparison and to demonstrate a clear advantage, we must consider the entire quantum phase, from state preparation (which introduces additional execution time) to the correction phase. In fact, we observe in Figure 12 that as the number of qubits increases, the execution time for state preparation inevitably increases.

\begin{figure}[H]
  \centering
    \includegraphics[width=1\textwidth]{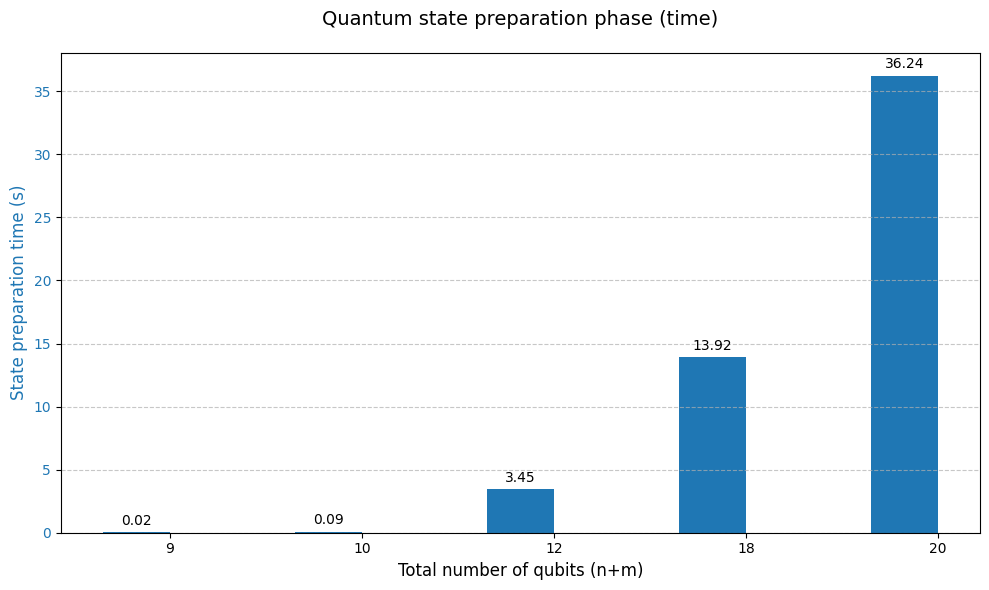}
    \caption{The presence of the state preparation and initialize the circuit cause additional time.}
\end{figure}

For this reason, if we aim to verify an actual temporal advantage, it is necessary to consider the entire quantum pipeline along with the correction phase, which introduces additional complexity. In Figure 13, the execution times of the correction procedure are plotted, and the results are equally satisfactory. Despite the faster growth, the quantum method's curve remains consistently below the classical one.
\begin{figure}[H]
  \centering
    \includegraphics[width=1\textwidth]{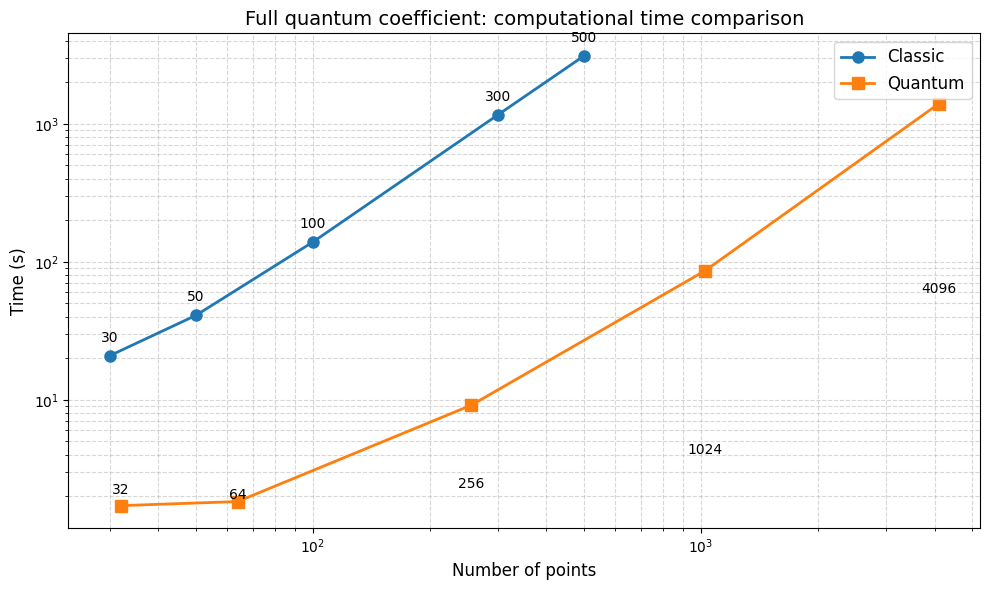}
    \caption{Runtime comparison between the quantum (with correction phase) and classical methods. The quantum approach exhibits faster growth but remains more efficient overall.}
\end{figure}

Regarding spatial complexity, the same observations and distinctions made for the quantum case apply. The results obtained from the QFT application phase alone are shown in Figure 13.

\begin{figure}[H]
  \centering
    \includegraphics[width=1\textwidth]{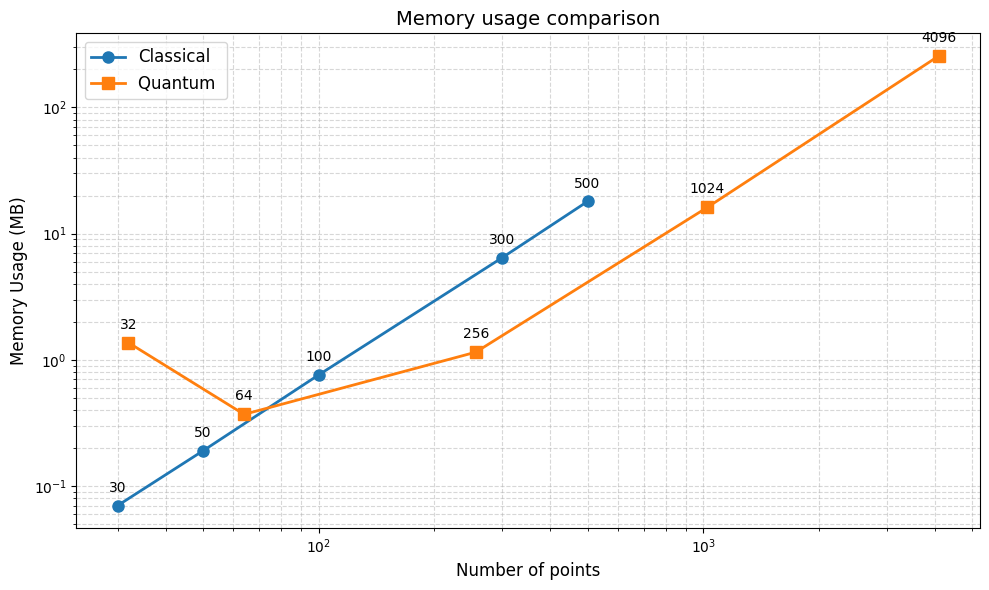}
    \caption{Memory usage is shown for the QFT step, excluding initialization and correction.}
\end{figure}
 For a fair comparison, an additional step defining the function must be included in the classical case, whereas the quantum approach already accounts both circuit initialization and coefficient correction.
Figure 15 reveals a slight memory advantage, though not as pronounced as in the temporal complexity case. More optimized spectral methods could further narrow the gap between the two curves. Nevertheless, the hybrid quantum approach does not require excessive memory usage, demonstrating a modest advantage.
\begin{figure}[H]
  \centering
    \includegraphics[width=1\textwidth]{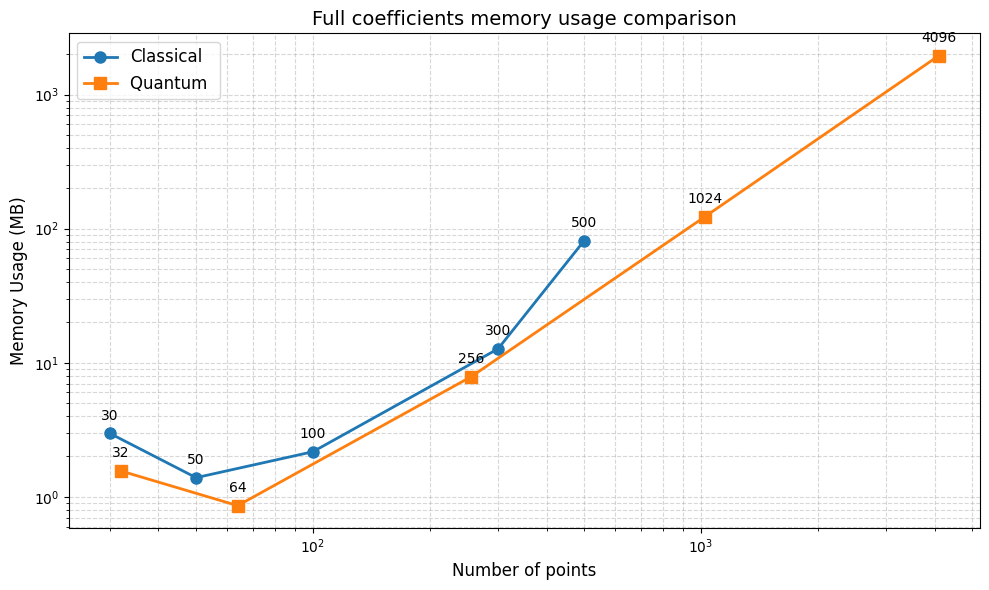}
    \caption{Comparative memory footprint of hybrid quantum vs. classical methods. Despite accounting for circuit initialization and coefficient correction, the quantum approach maintains a modest memory advantage over classical spectral methods.}
\end{figure}
In summary, the goal of performing the coefficient computation using quantum methods has proven to be successful, especially in terms of overall reduction in computational time.

\section{Conclusions}
The proposed method demonstrates a clear improvement, particularly in terms of computational time, for the estimation of the solution of the Poisson equation, especially when dealing with a large number of points. The objective of reducing computational time and memory usage, while providing accurate estimations for the solution of the Poisson equation through the computation of Fourier coefficients in a quantum framework, was convincingly achieved.

Nevertheless, future works could aim to compute the correction phase in quantum. However, it is important to acknowledge that the implementation does not rely solely on quantum computing; rather, it requires an interoperability framework with high-performance computing (HPC) systems. These systems would help accelerate the most computationally demanding benchmarking tasks of the algorithm, as discussed in the final section.

Furthermore, although the results obtained are promising, it is crucial to acknowledge that the experiments were conducted within a simulated environment.

In conclusion, our approach constitutes a significant contribution toward efficient solutions for partial differential equations, providing a solid foundation for tackling more complex problems and exploring novel interdisciplinary applications.

\nocite{*}
\bibliography{bibliography}

\begin{thebibliography}{29}
\providecommand{\natexlab}[1]{#1}
\providecommand{\url}[1]{\texttt{#1}}
\expandafter\ifx\csname urlstyle\endcsname\relax
  \providecommand{\doi}[1]{doi: #1}\else
  \providecommand{\doi}{doi: \begingroup \urlstyle{rm}\Url}\fi

\bibitem[Cochran et~al.(1967)Cochran, Cooley, Favin, Helms, Kaenel, Lang, Maling, Nelson, Rader, and Welch]{cochran1967fast}
William~T Cochran, James~W Cooley, David~L Favin, Howard~D Helms, Reginald~A Kaenel, William~W Lang, George~C Maling, David~E Nelson, Charles~M Rader, and Peter~D Welch.
\newblock What is the fast fourier transform?
\newblock \emph{Proceedings of the IEEE}, 55\penalty0 (10):\penalty0 1664--1674, 1967.

\bibitem[Daley et~al.(2022)Daley, Bloch, Kokail, Flannigan, Pearson, Troyer, and Zoller]{daley2022practical}
Andrew~J Daley, Immanuel Bloch, Christian Kokail, Stuart Flannigan, Natalie Pearson, Matthias Troyer, and Peter Zoller.
\newblock Practical quantum advantage in quantum simulation.
\newblock \emph{Nature}, 607\penalty0 (7920):\penalty0 667--676, 2022.

\bibitem[Davis and Rabinowitz(2007)]{davis2007methods}
Philip~J Davis and Philip Rabinowitz.
\newblock \emph{Methods of numerical integration}.
\newblock Courier Corporation, 2007.

\bibitem[Endo et~al.(2021)Endo, Cai, Benjamin, and Yuan]{endo2021hybrid}
Suguru Endo, Zhenyu Cai, Simon~C Benjamin, and Xiao Yuan.
\newblock Hybrid quantum-classical algorithms and quantum error mitigation.
\newblock \emph{Journal of the Physical Society of Japan}, 90\penalty0 (3):\penalty0 032001, 2021.

\bibitem[Evans(2010)]{evans10}
Lawrence~C. Evans.
\newblock \emph{Partial differential equations}.
\newblock American Mathematical Society, Providence, R.I., 2010.
\newblock ISBN 9780821849743 0821849743.

\bibitem[Faj et~al.(2023)Faj, Peng, Wahlgren, and Markidis]{faj2023quantum}
Jennifer Faj, Ivy Peng, Jacob Wahlgren, and Stefano Markidis.
\newblock Quantum computer simulations at warp speed: Assessing the impact of gpu acceleration: A case study with ibm qiskit aer, nvidia thrust \& cuquantum.
\newblock In \emph{2023 IEEE 19th International Conference on e-Science (e-Science)}, pages 1--10. IEEE, 2023.

\bibitem[Feng and Zhao(2020)]{feng2020fft}
Hongsong Feng and Shan Zhao.
\newblock Fft-based high order central difference schemes for three-dimensional poisson's equation with various types of boundary conditions.
\newblock \emph{Journal of Computational Physics}, 410:\penalty0 109391, 2020.

\bibitem[Freund and Wilson(2003)]{freund2003statistical}
Rudolf~J Freund and William~J Wilson.
\newblock \emph{Statistical methods}.
\newblock Elsevier, 2003.

\bibitem[Gautschi(2011)]{gautschi2011numerical}
Walter Gautschi.
\newblock \emph{Numerical analysis}.
\newblock Springer Science \& Business Media, 2011.

\bibitem[Gottlieb and Orszag(1977)]{gottlieb1977numerical}
David Gottlieb and Steven~A Orszag.
\newblock \emph{Numerical analysis of spectral methods: theory and applications}.
\newblock SIAM, 1977.

\bibitem[Javadi-Abhari et~al.(2024)Javadi-Abhari, Treinish, Krsulich, Wood, Lishman, Gacon, Martiel, Nation, Bishop, Cross, Johnson, and Gambetta]{javadiabhari2024quantumcomputingqiskit}
Ali Javadi-Abhari, Matthew Treinish, Kevin Krsulich, Christopher~J. Wood, Jake Lishman, Julien Gacon, Simon Martiel, Paul~D. Nation, Lev~S. Bishop, Andrew~W. Cross, Blake~R. Johnson, and Jay~M. Gambetta.
\newblock Quantum computing with qiskit, 2024.
\newblock URL \url{https://arxiv.org/abs/2405.08810}.

\bibitem[Kerr et~al.(2022)Kerr, Gonzalez-Parra, and Sherman]{article}
Gilbert Kerr, Gilberto Gonzalez-Parra, and Michelle Sherman.
\newblock A new method based on the laplace transform and fourier series for solving linear neutral delay differential equations.
\newblock \emph{Applied Mathematics and Computation}, 420:\penalty0 126914, 05 2022.
\newblock \doi{10.1016/j.amc.2021.126914}.

\bibitem[Kitaev et~al.(2002)Kitaev, Shen, and Vyalyi]{kitaev2002classical}
Alexei~Yu Kitaev, Alexander Shen, and Mikhail~N Vyalyi.
\newblock \emph{Classical and quantum computation}.
\newblock Number~47. American Mathematical Soc., 2002.

\bibitem[Krommer and Ueberhuber(1998)]{krommer1998_ComputationalIntegration_ueberhuber}
Arnold~R Krommer and Christoph~W Ueberhuber.
\newblock \emph{Computational integration}.
\newblock SIAM, 1998.

\bibitem[Lloyd(2003)]{lloyd2003hybrid}
Seth Lloyd.
\newblock Hybrid quantum computing.
\newblock In \emph{Quantum information with continuous variables}, pages 37--45. Springer, 2003.

\bibitem[Musk(2020)]{9198106}
D.~R. Musk.
\newblock A comparison of quantum and traditional fourier transform computations.
\newblock \emph{Computing in Science \& Engineering}, 22\penalty0 (6):\penalty0 103--110, 2020.
\newblock \doi{10.1109/MCSE.2020.3023979}.

\bibitem[Patan{\`e}(2018)]{patane2018laplacian}
Giuseppe Patan{\`e}.
\newblock Laplacian spectral basis functions.
\newblock \emph{Computer aided geometric design}, 65:\penalty0 31--47, 2018.

\bibitem[Patil and Prasad(2013)]{patil2013numerical}
Parag~V Patil and Dr~JSVR~Krishna Prasad.
\newblock Numerical solution for two dimensional laplace equation with dirichlet boundary conditions.
\newblock \emph{IOSR Journal of Mathematics}, 6\penalty0 (4):\penalty0 66--75, 2013.

\bibitem[Piessens et~al.(2012)Piessens, de~Doncker-Kapenga, {\"U}berhuber, and Kahaner]{piessens2012quadpack}
Robert Piessens, Elise de~Doncker-Kapenga, Christoph~W {\"U}berhuber, and David~K Kahaner.
\newblock \emph{Quadpack: a subroutine package for automatic integration}, volume~1.
\newblock Springer Science \& Business Media, 2012.

\bibitem[Sabetghadam et~al.(2009)Sabetghadam, Sharafatmandjoor, and Norouzi]{sabetghadam2009fourier}
Feriedoun Sabetghadam, Shervin Sharafatmandjoor, and Farhang Norouzi.
\newblock Fourier spectral embedded boundary solution of the poisson’s and laplace equations with dirichlet boundary conditions.
\newblock \emph{Journal of Computational Physics}, 228\penalty0 (1):\penalty0 55--74, 2009.

\bibitem[Scherer(2019)]{scherer2019mathematics}
Wolfgang Scherer.
\newblock \emph{Mathematics of quantum computing}, volume~11.
\newblock Springer, 2019.

\bibitem[Shen et~al.(2011)Shen, Tang, and Wang]{shen2011spectral}
Jie Shen, Tao Tang, and Li-Lian Wang.
\newblock Spectral methods: algorithms, analysis and applications.
\newblock 41, 2011.

\bibitem[Shortley and Weller(1938)]{shortley1938numerical}
George~H Shortley and Royal Weller.
\newblock The numerical solution of laplace's equation.
\newblock \emph{Journal of Applied Physics}, 9\penalty0 (5):\penalty0 334--348, 1938.

\bibitem[Steane(1998)]{steane1998quantum}
Andrew Steane.
\newblock Quantum computing.
\newblock \emph{Reports on Progress in Physics}, 61\penalty0 (2):\penalty0 117, 1998.

\bibitem[Strauss(2007)]{strauss2007partial}
W.A. Strauss.
\newblock \emph{Partial Differential Equations: An Introduction}.
\newblock Wiley, 2007.
\newblock ISBN 9780470054567.
\newblock URL \url{https://books.google.it/books?id=m2hvDwAAQBAJ}.

\bibitem[Thom{\'e}e(2001)]{thomee2001finite}
Vidar Thom{\'e}e.
\newblock From finite differences to finite elements a short history of numerical analysis of partial differential equations.
\newblock In \emph{Numerical analysis: Historical developments in the 20th century}, pages 361--414. Elsevier, 2001.

\bibitem[Weinstein et~al.(2001)Weinstein, Pravia, Fortunato, Lloyd, and Cory]{PhysRevLett.86.1889}
Y.~S. Weinstein, M.~A. Pravia, E.~M. Fortunato, S.~Lloyd, and D.~G. Cory.
\newblock Implementation of the quantum fourier transform.
\newblock \emph{Phys. Rev. Lett.}, 86:\penalty0 1889--1891, Feb 2001.
\newblock \doi{10.1103/PhysRevLett.86.1889}.
\newblock URL \url{https://link.aps.org/doi/10.1103/PhysRevLett.86.1889}.

\bibitem[Widder(1976)]{widder1976heat}
David~Vernon Widder.
\newblock \emph{The heat equation}, volume~67.
\newblock Academic Press, 1976.

\bibitem[Williams(2011)]{Williams2011}
Colin~P. Williams.
\newblock \emph{Quantum Gates}, pages 51--122.
\newblock Springer London, London, 2011.
\newblock ISBN 978-1-84628-887-6.
\newblock \doi{10.1007/978-1-84628-887-6_2}.
\newblock URL \url{https://doi.org/10.1007/978-1-84628-887-6_2}.

\end{thebibliography}

\end{document}